\documentclass[reqno,tbtags,intlimits,a4paper,oneside,12pt]{amsart}
\usepackage[cp1251]{inputenc}%
\usepackage[english]{babel}
\usepackage{amssymb,upref,calc,url}
\usepackage[mathscr]{eucal}
\usepackage{graphicx}
\usepackage{amsmath,amscd,amsthm,verbatim}
\usepackage[all]{xy}
\usepackage[in]{fullpage}
\newtheorem{thm}{\hskip\parindent Theorem}[section]

\newtheorem{cor}[thm]{\hskip\parindent Corollary}

\theoremstyle{definition}

\DeclareMathOperator{\Jac}{Jac}

\begin{document}

\title{Generators in the field of hyperelliptic functions}
\author{E.\,Yu.~Bunkova}
\address{Steklov Mathematical Institute of Russian Academy of Sciences, Moscow, Russia}
\email{bunkova@mi-ras.ru}

\begin{abstract}
We consider the field of hyperelliptic functions defined for a family of hyperelliptic curves as rational functions in some special functions from Kleinian functions theory. We compare our definition with the classical one. We provide details and references for the result that the field of hyperelliptic functions for a family of hyperelliptic curves of genus $g$ is isomorphic to the field of rational functions with $3g$ generators. The main result of the present work is that there are no algebraic relations between these generators.
\end{abstract} \maketitle

\section{Introduction} \label{S0}

Let $g \in \mathbb{N}$. The space of hyperelliptic functions of genus $g$ is constructed starting from a curve with the affine part of the form
\begin{equation} \label{V1}
\mathcal{V}_\lambda = \{(x, y)\in\mathbb{C}^2:
y^2 = x^{2g+1} + \lambda_4 x^{2 g - 1} + \lambda_6 x^{2 g - 2} + \cdots + \lambda_{4 g} x + \lambda_{4 g + 2}\}.
\end{equation}
We denote the parameters of the equation in \eqref{V1} by $\lambda = (\lambda_4, \lambda_6,\dots, \lambda_{4 g}, \lambda_{4 g + 2}) \in \mathbb{C}^{2g}$.
A hyperelliptic curve is defined by a specialization of the parameters $\lambda$ such that the equation 
\[
 x^{2g+1} + \lambda_4 x^{2 g - 1} + \lambda_6 x^{2 g - 2} + \cdots + \lambda_{4 g} x + \lambda_{4 g + 2} = 0
\] 
has no multiple roots. We denote by $\mathcal{B}_g$ the subspace of parameters $\mathcal{B}_g \subset \mathbb{C}^{2g}$ of such curves.

In classical works hyperelliptic functions of genus $g$ are constructed starting from a curve~\eqref{V1} with some parameters $\lambda \in \mathcal{B}_g$, considering the set of parameters $\lambda$ of the curve~$\mathcal{V}_\lambda$ as constants. In the present work, following the idea of \cite{BL}, we consider the whole family of curves of the form \eqref{V1} with all the parameters $\lambda \in \mathbb{C}^{2g}$ simultaniously to explore the dependence of hyperelliptic functions of genus $g$ on the parameters $\lambda$. We give the relevant definitions in \S \ref{S1} based on the theory of hyperelliptic Kleinian functions, see~\cite{K88},~\cite{Baker},~\cite{BEL-97},~\cite{BEL},~\cite{BEL-12},~\cite{BEL18},~and~\cite{BMulti}.  Note that in the classical case for fixed parameters $\lambda \in \mathcal{B}_g$ these definitions give the classical definitions of hyperelliptic functions, see \cite{BL} and also \cite{WW} for elliptic functions.
We denote $\mathcal{B}_g = \mathbb{C}^{2g} \setminus\Sigma_g$, where $\Sigma_g$ is the discriminant hypersurface of the family of hyperelliptic curves \eqref{V1}.

In \S \ref{S2} we describe some relations between hyperelliptic functions based on the works~\cite{BEL-97} and \cite{BS}.
In \S \ref{S3} we use these relations and results from our works \cite{BB24}, \cite{BB23}, \cite{BB22}, \cite{B3} to show that there are no algebraic relations between the $3g$ functions that we choose as generators in the field of hyperelliptic functions.

The form \eqref{V1} for a family of hyperelliptic curves and the corresponding field of functions were studied in the work \cite{BL}. Moreover, in this work a more general case of $(n,s)$-curves was considered. Here we limit ourselves to the case of the form \eqref{V1}, that is the case of~$(2,2g+1)$-curves in the notation of the work \cite{BL}.

In this work all the equations will be homogeneous with respect to some grading of the variables. Generally the indices of the variables will represent their grading. We denote the grading of a variable $x$ by $|x|$. For instance in \eqref{V1} we have $|x| = 2$, $|y| = 2g+1$, and~$|\lambda_s| = s$ where $s \in \{4, 6, \ldots, 4 g, 4 g + 2\}$.

\section{The field of hyperelliptic functions} \label{S1}

We consider the complex space $\mathbb{C}^g$ and we denote the coordinates in this space by $z_s$ for odd integers $s$ from $1$ to $2g-1$, where $|z_s| = s$. We denote $z = (z_1, z_3,\dots, z_{2g-1}) \in \mathbb{C}^g$.  The vector $\omega \in \mathbb{C}^g$ is called a period of a meromorphic function $f$ on $\mathbb{C}^g$ if $f(z+\omega) = f(z)$ for all~$z \in \mathbb{C}^g$. If a meromorphic function $f$ has $2g$ independent periods in~$\mathbb{C}^g$, then it is called an Abelian function. Thus, an Abelian function is a meromorphic function on the complex torus $\mathbb{C}^g\!/\Gamma$, where $\Gamma$ is a lattice formed by the periods. In this work we are interested in Abelian functions in $\mathbb{C}^g$ such that their lattices of periods are determined by a hyperelliptic curve of the form \eqref{V1}, see \cite{BL} for the definitions and \cite{BEL} or \S 2 in \cite{BB24} for the explicit construction of such a lattice. In this case the torus $\mathbb{C}^g\!/\Gamma$ is called the Jacobi manifold (or, shortly, Jacobian) of the curve. The universal bundle of Jacobians of hyperelliptic curves is the bundle $\pi\colon \mathcal{U}_g \to \mathcal{B}_g$ with fiber over $\lambda \in \mathcal{B}_g$ the Jacobian of the curve~$\mathcal{V}_\lambda$. Hyperelliptic functions of genus~$g$ are  functions on the total space $\mathcal{U}_g$ of the bundle $\pi$. A restriction of such a function for each $\lambda \in \mathcal{B}_g$ is an Abelian function with lattice of periods determined by the curve~$\mathcal{V}_\lambda$. 

We define the field of hyperelliptic functions based on the theory of hyperelliptic Kleinian functions, see 
\cite{K88}, \cite{Baker}, \cite{BEL-97}, \cite{BEL}, \cite{BEL-12},~\cite{BEL18}, and~\cite{BMulti}. Let 
$\sigma(z, \lambda)$ be the genus $g$ hyperelliptic sigma function. This function depends on the variables $z = (z_1, z_3,\dots, z_{2g-1}) \in \mathbb{C}^g$ and the parameters $\lambda = (\lambda_4, \lambda_6,\dots, \lambda_{4 g}, \lambda_{4 g + 2}) \in \mathbb{C}^{2g}$ that correspond to the parameters of the curve \eqref{V1}. Set
\begin{equation}
\wp_{k_1,\dots, k_n} = - \partial_{k_1} \cdots \partial_{k_n} \ln \sigma(z, \lambda), \label{V2}
\end{equation}
where $n \geqslant 2$ and $k_s \in \{ 1, 3,\dots, 2 g - 1\}$. From the theory of hyperelliptic Kleinian functions, see Chapter 3 in \cite{BEL} and equation (6.3) in \cite{BB24}, it follows that the functions $\wp_{k_1,\dots, k_n}$ satisfy the properties needed for being a hyperelliptic function, namely a restriction of such a function for each $\lambda \in \mathcal{B}_g$ is an Abelian function with lattice of periods determined by the curve~$\mathcal{V}_\lambda$. Naturally, rational functions in any functions $\wp_{k_1,\dots, k_n}$ with coefficients in $\mathbb{C}[\lambda]$ satisfy these properties as well. Following the works \cite{BB24} and \cite{BB23}, we define hyperelliptic functions as rational functions in $\wp_{k_1,\dots, k_n}$ with coefficients in $\mathbb{C}[\lambda]$ and denote the field of hyperelliptic functions by $\mathcal{F}$. This definition is motivated by the fact that for any $\lambda \in \mathcal{B}_g$ any Abelian function on the Jacobian of the curve~$\mathcal{V}_\lambda$ is equal to the restriction of a function from $\mathcal{F}$ to this Jacobian, see property (3) in Section 1.1 in \cite{BL}. By Corollary 6.2 in \cite{BB24}, the field $\mathcal{F}$ is isomorphic to the field of rational functions over~$\mathbb{Q}$ in $3g$ generators $\wp_{1,k}, \wp_{1,1,k}, \wp_{1,1,1,k}$, where $k \in \{1,3, \ldots, 2g-1 \}$. We~give the details in~\S \ref{S3} of the present work.

In the definition of the field of hyperelliptic functions $\mathcal{F}$ there are no additional conditions on the parameters $\lambda$, so one can consider such functions both with parameters~$\lambda \in \mathcal{B}_g$ and with parameters $\lambda \in \Sigma_g$. For example, we have $0 \in \Sigma_g$ for any genus $g \in \mathbb{N}$. We give details on the construction of the function $\sigma(z, \lambda)$ in \eqref{V2} in the special case $\lambda = 0$ in~\cite{BB21}. Note that for $\lambda \in \Sigma_g$ a function $f \in \mathcal{F}$ restricted to $\mathbb{C}^g \times \lambda$ is not necessarily an Abelian function, therefore it is not a hyperelliptic function in the classical sense.

Let us compare this definition with an older definition of the field $\mathcal{F}$, see for instance~\S 4 in \cite{BB22}, \S 1 in \cite{B3}, \S 1 in \cite{BB21}, or Section 3 in the review~\cite{BM}, which defines the field $\mathcal{F}$ as the field of functions on $\mathbb{C}^g \times \mathcal{B}_g$ such that a restriction of $f \in \mathcal{F}$ to $\mathbb{C}^g \times \lambda$ is an Abelian function on $\mathbb{C}^g$ with periods defined by the curve $\mathcal{V}_\lambda$. The review states that such a field is isomorphic to the field of rational functions in $\wp_{k_1,\dots, k_n}$. However, counterexamples to this statement, such as the function~$f = \exp(\lambda_4)$, can be given.

\section{Relations in the field of hyperelliptic functions} \label{S2}

The relations described in this section originate from the explicit algebraic relations between the generating functions in the field of meromorphic functions on Jacobians of genus $g$ hyperelliptic curves, see \S 3 in \cite{BEL-97}. Here we give the formulas in the special case of the curve $\mathcal{V}_\lambda$ defined by \eqref{V1}. We use the modern notations that are consistent throughout the work. We use the notation $\delta_{\cdot,\cdot}$ for the Kronecker delta function. We set $\wp_{i,j} = 0$ for~$i \geqslant 2g+1$ or $j \geqslant 2g+1$.

By Corollary~3.1.2 in~\cite{BEL-97} for any $i\in \{1, 3,\dots, 2 g - 1\}$ we have the relation:
\begin{equation}
\wp_{1,1,1,i} = 6 \wp_{1,1} \wp_{1,i} - 2 \wp_{3, i} + 6 \wp_{1, i+2} + 2 \lambda_{4} \delta_{i,1}. \label{BEL1}
\end{equation}

By Theorem 3.2 in \cite{BEL-97} for any $i, j \in \{1, 3,\dots, 2 g - 1\}$ we have the relation:
\begin{multline}
\wp_{1,1,i} \wp_{1,1,j} = 4 \wp_{1,1} \wp_{1,i} \wp_{1,j} + 4 \wp_{1,i+2} \wp_{1,j} + 4 \wp_{1,i} \wp_{1, j+2} - 2 \wp_{3,i} \wp_{1,j} - 2 \wp_{1,i} \wp_{3,j} - \\ - 2 (\wp_{i+4, j} - 2 \wp_{i+2, j+2} + \wp_{i, j+4})
 + 2 \lambda_4 (\delta_{i,1} \wp_{1,j} + \wp_{1,i} \delta_{j,1}) + 2 \lambda_{i+j+4} (2 \delta_{i,j} + \delta_{i-2,j} + \delta_{i, j-2}). \label{BEL2}
\end{multline}

The work \cite{BEL-97} considers the case of a hyperelliptic curve $\mathcal{V}_\lambda$ with fixed parameters $\lambda \in \mathcal{B}_g$. In this case
the relations \eqref{BEL2} are considered as relations between generators in the field of meromorphic functions on the Jacobian of $\mathcal{V}_\lambda$, namely 
between the functions $\wp_{i,j}$, where  $i, j \in \{1, 3,\dots, 2 g - 1\}$, and $\wp_{1,1,k}$, where $k \in \{1, 3,\dots, 2 g - 1\}$. This set of functions is considered as generating functions, because for fixed parameters $\lambda \in \mathcal{B}_g$ any meromorphic function on the Jacobian of the curve $\mathcal{V}_\lambda$ can be presented as a rational function in these generators. The parameters $\lambda \in \mathcal{B}_g$ of the curve are considered as parameters for the relations \eqref{BEL2}.

In the case of hyperelliptic functions depending on the parameters $\lambda$ of the family of hyperelliptic curves~$\mathcal{V}_\lambda$, see definition in \S \ref{S1}, we consider the parameters $\lambda$ as functions in $\mathcal{F}$. We use the relations~\eqref{BEL2} along with the relations \eqref{BEL1} to express the functions $\lambda_s$ for $s \in \{4, 6, \ldots, 4g, 4g+2\}$ and~$\wp_{i,j}$, where  $i, j \in \{3, 5,\dots, 2 g - 1\}$, as functions in the $3g$ generators $\wp_{1,k}, \wp_{1,1,k}, \wp_{1,1,1,k}$ for $k \in \{1, 3,\dots, 2 g - 1\}$ of the field of hyperelliptic functions, see \S \ref{S3}. Note that in contrast with the case of the work \cite{BEL-97}, there are no relations between the generators $\wp_{1,k}, \wp_{1,1,k}, \wp_{1,1,1,k}$ for $k \in \{1, 3,\dots, 2 g - 1\}$ in the field of hyperelliptic functions, see Theorem \ref{pt} in \S \ref{S3} of the present work.

In the work \cite{BS} in Theorem 4.1 there is a convenient expression for the parameters~$\lambda$ as functions in $\mathcal{F}$, see Corollary 5.2 and Lemma 5.3 in \cite{BB22} for modern notations and a~detailed proof. Namely, we have the relation
\begin{equation} \label{L1}
 4 {\bf m}(\xi) = {\bf b}_2(\xi)^2 + 2 {\bf b}_3(\xi) (1 - {\bf b}_1(\xi)) + 4 (\xi^{-1} + 2 \wp_{1,1}) (1 - {\bf b}_1(\xi))^2, 
\end{equation}
where 
\begin{align*}
{\bf b}_1(\xi) &= \sum_{i=1}^g \wp_{1,2i-1} \xi^i, & 
{\bf b}_2(\xi) &= \sum_{i=1}^g \wp_{1,1,2i-1} \xi^i, &
{\bf b}_3(\xi) &= \sum_{i=1}^g \wp_{1,1,1,2i-1} \xi^i,
\end{align*}
\begin{align*}
{\bf m}(\xi) &= \xi^{-1} + \sum_{i=1}^{2g} \lambda_{2i+2} \xi^i,
\end{align*}
and $\xi$ is a free parameter, so the equation \eqref{L1} gives relations expressing the parameters~$\lambda_s$ for $s \in \{4, 6, \ldots, 4g, 4g+2\}$ in its left hand side as polynomials in the $3g$ generators $\wp_{1,k}, \wp_{1,1,k}, \wp_{1,1,1,k}$ for $k \in \{1, 3,\dots, 2 g - 1\}$ in the field of hyperelliptic functions. Other expressions with the same result can be found in Lemma 6.3 in \cite{BB23} and Corollary 5.5 in~\cite{B3}.

\section{Generators in the field of hyperelliptic functions} \label{S3}

Theorem 4.1 in \cite{BB23} states that all the functions~\eqref{V2} and all the parameters $\lambda_s$ of the curve~\eqref{V1} can be expressed as polynomials with coefficients in $\mathbb{Q}$ of the $3g$ generators $\wp_{1,k}, \wp_{1,1,k}, \wp_{1,1,1,k}$ for $k \in \{1, 3,\dots, 2 g - 1\}$. The proof of this Theorem is obtained from the relations~\eqref{BEL1}, \eqref{BEL2} and \eqref{L1}. The relation \eqref{L1} gives expressions for the parameters~$\lambda_s$ in our $3g$ generators. Using these expressions, we get expressions for $\wp_{3,k}$ for~$k \in \{3,5,\dots, 2 g - 1\}$ from \eqref{BEL1} and for $\wp_{i,j}$ for $i,j \in \{5,7,\dots, 2 g - 1\}$ from \eqref{BEL2}. All the other functions in~\eqref{V2} either belong to our set of generators or can be obtained from relations that are in turn obtained by differentiating the relations~\eqref{BEL1} and \eqref{BEL2} with respect to the variables $z_i$.

A direct Corollary of Theorem 4.1 in \cite{BB23}, obtained by using the definition of the field~$\mathcal{F}$ from~\S \ref{S1}, is that any function $f \in \mathcal{F}$ can be expressed as a rational function with coefficients in $\mathbb{Q}$ of the $3g$ generators $\wp_{1,k}, \wp_{1,1,k}, \wp_{1,1,1,k}$ for $k \in \{1, 3,\dots, 2 g - 1\}$.

\begin{cor}[Corollary 6.2 in \cite{BB24}] \label{c62}
The field $\mathcal{F}$ is isomorphic to the field of rational functions over~$\mathbb{Q}$ in $3g$ generators $\wp_{1,k}, \wp_{1,1,k}, \wp_{1,1,1,k}$, where $k \in \{1,3, \ldots, 2g-1 \}$.
\end{cor}

In the following theorem we consider the relations between these generators.

\begin{thm} \label{pt}
There is no algebraic relations in the field $\mathcal{F}$ of hyperelliptic functions of genus $g$ between any of the $3g$ functions $\wp_{1,k}$, $\wp_{1,1,k}$, $\wp_{1,1,1,k}$ for $k \in \{1, 3,\dots, 2 g - 1\}$, considered as functions in the total space $\mathcal{U}_g$ of the universal bundle of Jacobians of genus~$g$ hyperelliptic curves.
\end{thm}

The proof follows from Theorem 5.3 in \cite{B3} and Corollary 3.2.1 in \cite{BEL-97}.

To explain the notation of Theorem 5.3 in \cite{B3} we give the statement of Corollary~5.2 in~\cite{B3} in a slightly modified notation to be consistent with the definitions above. Corollary~5.2 in \cite{B3} considers the map $\psi: \mathcal{U}_g \dashrightarrow \mathbb{C}^{ \frac{g (g+9)}{2}}$, where in $\mathbb{C}^{ \frac{g (g+9)}{2}}$ we denote
 the coordinates by
 $(b, w, \lambda)$. Here for $b = (b_1, b_2, b_3)= (b_{1,j}, b_{2,j}, b_{3,j}) \in \mathbb{C}^{3g}$ we have  \mbox{$j \in \{1, 3, \ldots, 2 g -1\}$,}
 for $w = (w_{k,l}) \in \mathbb{C}^{\frac{g (g-1)}{2}}$ we have $k,l \in \{3, 5, \ldots, 2 g -1\}$, $k \leqslant l$,
 and for $\lambda = (\lambda_s) \in \mathbb{C}^{2 g}$ we have $s \in \{4, 6, \ldots, 4 g, 4 g + 2\}$.
 We set 
\begin{equation} \label{EP}
\psi: (z, \lambda) \mapsto (b_{1,j}, b_{2,j}, b_{3,j}, w_{k,l}, \lambda_s) = (\wp_{1,j}, \wp_{1,1,j}, \wp_{1,1,1,j}, \wp_{k,l}, \lambda_s).
\end{equation}
 
Corollary~5.2 in \cite{B3} states that the image of $\psi$ lies in $\mathcal{S}_g \subset \mathbb{C}^{ \frac{g (g+9)}{2}}$, where $\mathcal{S}_g$ is determined by the set of $\frac{g (g+3)}{2}$ equations that are the images of the equations \eqref{BEL1} and~\eqref{BEL2} under the map $\psi$, these equations are given explicitly in~\cite{B3}.

\begin{thm}[Theorem 5.3 in \cite{B3}] \label{thm3}
The projection $\pi_1\colon \mathbb{C}^{\frac{g (g+9)}{2}} \to \mathbb{C}^{3 g}$ on the first $3 g$ coordinates gives the isomorphism $\mathcal{S}_g \simeq \mathbb{C}^{3g}$.
Therefore, the coordinates $b$ uniformize $\mathcal{S}_g$.
\end{thm}

We revisit this Theorem in Theorem 6.2 in \cite{BB22}.
The work \cite{B3} illustrates this Theorem by a commutative diagram
\begin{equation} \label{D1}
\xymatrix{
	 & \mathbb{C}^{\frac{g (g+9)}{2}} \ar@{<-_{)}}[d] \ar@{=}[r] & \mathbb{C}^{3 g} \times \mathbb{C}^{\frac{g (g-1)}{2}} \times \mathbb{C}^{2 g} \ar[d]^{\pi_1}  \\
	\mathcal{U}_g \ar[d]^{\pi} \ar@{-->}[r]^{\psi} & \mathcal{S}_g \ar@{=}[r]^{\sim} & \mathbb{C}^{3 g} \ar[d]^{p}\\
	\mathcal{B}_g \ar@{^{(}->}[rr] & & \mathbb{C}^{2g}\\
	},
\end{equation}
where we can obtain the polynomial map $p$ from the relation \eqref{L1} and all the complex spaces have the coordinates according to their dimensions as described above. The Theorem considers the upper half of diagram \eqref{D1}.

We recall that $\mathcal{B}_g = \mathbb{C}^{2g} \setminus\Sigma_g$, where $\Sigma_g$ is the discriminant hypersurface of the family of hyperelliptic curves \eqref{V1}. The commutativity of the lower half of diagram \eqref{D1} implies that the image of $\psi$ does not intersect with the preimage of $\Sigma_g \subset \mathbb{C}^{2g}$ under the polynomial map $p$. Therefore the image of $\psi$ in $\mathcal{S}_g \simeq \mathbb{C}^{3g}$ lies in $\mathbb{C}^{3g} \backslash p^{-1}(\Sigma_g)$.

We recall that in the bundle $\pi\colon \mathcal{U}_g \to \mathcal{B}_g$ the fiber over $\lambda \in \mathcal{B}_g$ is the Jacobian of the curve~$\mathcal{V}_\lambda$. Following \cite{BEL-97} we denote this Jacobian by $\Jac(\mathcal{V}_\lambda)$ and by $(\sigma)$ we denote the divisor of zeroes of the function $\sigma(z,\lambda)$ in $\Jac(\mathcal{V}_\lambda)$. We consider the complex space~$\mathbb{C}^{\frac{g (g+3)}{2}}$ where, following the notation above, we denote the coordinates by $(b_1, b_2, w)$.
Here for~$(b_1, b_2) = (b_{1,j}, b_{2,j}) \in \mathbb{C}^{2g}$ we have \mbox{$j \in \{1, 3, \ldots, 2 g -1\}$,}
 for $w = (w_{k,l}) \in \mathbb{C}^{\frac{g (g-1)}{2}}$ we have $k,l \in \{3, 5, \ldots, 2 g -1\}$, $k \leqslant l$.
By Corollary 3.2.1 in \cite{BEL-97}, the map
 \[
  \varphi: \Jac(\mathcal{V}_\lambda) \backslash (\sigma) \to \mathbb{C}^{{\frac{g (g+3)}{2}}}, \qquad (z, \lambda) \mapsto (b_{1,j}, b_{2,j}, w_{k,l}) = (\wp_{1,j}, \wp_{1,1,j}, \wp_{k,l})
 \]
is a meromorphic embedding.

We consider the bundle $\pi_0 \colon \mathcal{U}_g^\circ \to \mathcal{B}_g$ associated with the universal bundle of Jacobians of hyperelliptic curves $\pi$. The fiber of the bundle $\pi_0$ over $\lambda \in \mathcal{B}_g$ is~$\Jac(\mathcal{V}_\lambda) \backslash (\sigma)$ and~$\mathcal{U}_g^\circ$ denotes the total space of the bundle $\pi_0$. Note that the map $\psi$ given by \eqref{EP} is well defined on $\mathcal{U}_g^\circ$ because in equation \eqref{V2} we see that the poles of the meromorphic functions involved in $\psi$ belong to the divisor of zeroes of the function $\sigma(z,\lambda)$ in $\Jac(\mathcal{V}_\lambda)$.

Set $\phi = \pi_1 \circ \psi$. As $\pi_1$ is the projection on the first three coordinates, from \eqref{EP} we get
\begin{equation} \label{EH}
\phi: (z, \lambda) \mapsto (b_{1,j}, b_{2,j}, b_{3,j}) = (\wp_{1,j}, \wp_{1,1,j}, \wp_{1,1,1,j}).
\end{equation}

Now to give the proof of Theorem \ref{pt} we consider the diagram
\begin{equation} \label{D2}
\xymatrix{
	\mathcal{U}_g^\circ \ar[d]^{\pi} \ar[r]^(.3){\phi} & \mathbb{C}^{3g} \backslash p^{-1}(\Sigma_g) \ar@{^{(}->}[r] \ar[d]^{p} & \mathbb{C}^{3 g} \ar[d]^{p}\\
	\mathcal{B}_g \ar@{=}[r] & \mathcal{B}_g \ar@{^{(}->}[r] & \mathbb{C}^{2g}\\
	}.
\end{equation}
For fixed $\lambda \in \mathcal{B}_g$, by Corollary 3.2.1 in \cite{BEL-97}, the map $\varphi$, and therefore the restriction of the map $\phi$ to the fiber over $\lambda \in \mathcal{B}_g$, is a meromorphic embedding. Therefore the map~$\phi$ gives a fibervise meromorphic embedding of $\mathcal{U}_g^\circ$ into $\mathbb{C}^{3g} \backslash p^{-1}(\Sigma_g) \subset \mathbb{C}^{3 g}$. Note that by Dubrovin--Novikov theorem, see \cite{DN}, the space $\mathcal{U}_g$ is birationally equivalent to~$\mathbb{C}^{3g}$, and therefore it is $3g$-dimensional. A more direct way to see this is that $\Jac(\mathcal{V}_\lambda)$ is $g$-dimensional as a factor of $\mathbb{C}^g$ by a lattice, while $\mathcal{B}_g$ is a $2g$-dimensional subset of $\mathbb{C}^{2g}$. In~the embedding $\phi$ of~$\mathcal{U}_g^\circ$ into $\mathbb{C}^{3 g}$ the coordinates $b$ in $\mathbb{C}^{3 g}$ are mapped into the $3g$ functions $\wp_{1,k}, \wp_{1,1,k}, \wp_{1,1,1,k}$, where $k \in \{1,3, \ldots, 2g-1 \}$, therefore there are no algebraic relations between these functions as functions in $\mathcal{U}_g^\circ$. Theorem \ref{pt} is proved.

The work \cite{B3} essentially uses this result while solving the problem of differentiation of hyperelliptic functions. In \S 6 of \cite{B3} the problem of finding the derivatives of the field $\mathcal{F}$ is replaced by the problem of finding polynomial vector fields in $\mathbb{C}^{3 g}$. In the associated construction the $3g$ functions $\wp_{1,k}, \wp_{1,1,k}, \wp_{1,1,1,k}$, where $k \in \{1,3, \ldots, 2g-1 \}$, are taken as generators of the field $\mathcal{F}$ and their algebraic independence is supposed. We obtain the following results:

\begin{thm} \label{TM}
The field $\mathcal{F}$ of hyperelliptic functions for a family of hyperelliptic curves~\eqref{V1} of genus $g$ is isomorphic to the field of rational functions in $\mathbb{C}^{3g}$.
\end{thm}

Let us consider the left part of the diagram \eqref{D2}
\begin{equation} \label{D3}
\xymatrix{
	\mathcal{U}_g^\circ \ar[d]^{\pi} \ar[r]^(.3){\phi} & \mathbb{C}^{3g} \backslash p^{-1}(\Sigma_g)  \ar[d]^{p} \\
	\mathcal{B}_g \ar@{=}[r] & \mathcal{B}_g \\
	}.
\end{equation}

\begin{thm} \label{TN}
The map~$\phi$ defined by \eqref{EH} gives a fibervise meromorphic embedding in diagram \eqref{D3}. For each $\lambda \in \mathcal{B}_g$ it maps $\Jac(\mathcal{V}_\lambda) \backslash (\sigma)$ into $p^{-1}(\lambda)$.
\end{thm}

\end{document}